%% file: uwVarNet_arxiv.tex
\documentclass{amsart}

\newtheorem{thm}{Theorem}[section]
\newtheorem{Proposition}[thm]{Proposition}

\newcounter{example}
\theoremstyle{definition}
\newtheorem{Example}[thm]{Example}

\numberwithin{equation}{section}

\usepackage{booktabs}
\usepackage{cite}

\usepackage[linesnumbered,ruled,vlined]{algorithm2e}
\usepackage{csquotes}
\newcommand{\Vvert}{{\vert\kern-0.25ex\vert\kern-0.25ex\vert}}
\usepackage{tikz,pgfplots}
\pgfplotsset{select coords between index/.style 2 args={
    x filter/.code={
        \ifnum\coordindex<#1\fi
        \ifnum\coordindex>#2\fi
    }
}}

\usepackage{textcomp}
\usepackage[textsize=tiny]{todonotes}

\newcommand\cB{\mathcal{B}}
\newcommand\cP{\mathcal{P}}

\newcommand\R{\mathbb{R}}
\newcommand\X{\mathbb{X}}

\begin{document}

\title[Wavelet-based PINNs]{A certified wavelet-based physics-informed neural network for the solution of parameterized partial differential equations}

\author{Lewin Ernst}
\address{Institute for Numerical Mathematics, Ulm University, Helmholtzstr. 20, 89081 Ulm, Germany}
\email{lewin.ernst@uni-ulm.de}

\author{Karsten Urban}
\address{Institute for Numerical Mathematics, Ulm University, Helmholtzstr. 20, 89081 Ulm, Germany}
\email{karsten.urban@uni-ulm.de}

\begin{abstract}
Physics Informed Neural Networks (PINNs) have frequently been used for the numerical approximation of Partial Differential Equations (PDEs). The goal of this paper is to construct PINNs along with a computable upper bound of the error, which is particularly relevant for model reduction of Parameterized PDEs (PPDEs). To this end, we suggest to use a weighted sum of expansion coefficients of the residual in terms of an adaptive wavelet expansion both for the loss function and an error bound. \\
This approach is shown here for elliptic PPDEs using both the standard variational and an optimally stable ultra-weak formulation. Numerical examples show a very good quantitative effectivity of the wavelet-based error bound.
\end{abstract}



\maketitle

\input{paper_uwVarNet.tex}

\bibliographystyle{ieeetr}
\bibliography{uwVarNet.bib}

\end{document}

%% file: paper_uwVarNet.tex

\section{Introduction}

Physics Informed Neural Networks (PINNs) have recently been introduced for the numerical solution of Partial Differential Equations (PDEs), see e.g.\ \cite{Raissi2019, Berner2021, Kharazmi2019, Khodayi-Mehr2020, Kharazmi2021, Jagtap2020, pang2019fpinns, cai2021physics}. Even though the range of PDEs that can be approximated seemingly well by PINNs is quite impressive, a rigorous a posteriori error control is at least not straightforward. Some attempts to bound the error and thus to certify PINN computations have been reported in the literature, e.g.\ \cite{BCP2,BCP,Shin2020v1, Shin2020v2, de2021error, mishra2022estimates1, mishra2022estimates2, Jiao2022}. Those PINNs use the residual of the PDE in the definition of the loss function to be minimized during the training phase. 

Several PINNs use the classical (i.e., point-wise) formulation of the PDE and produce point-wise approximations of the solution of a given PDE. In those cases, also the residual (used for the training) is evaluated point-wise, so that point-evaluation needs to be well-defined at least at those points. However, if the residual is required to be at least piecewise continuous, the solution of the PDE is implicitly required to have high regularity. Hence, PINNs requiring point-wise residuals work by construction on the classical form of the PDE.

However, defining the loss function in terms of the residual and hoping to obtain a certified PINN requires that the residual allows to define an upper bound for the error. This, in turn, requires that the formulation of the PDE used within the PINN is well-posed (the corresponding PDE operator needs to be an isomorphism), which is, at least in general, not true for the classical formulation of PDEs, \cite{Gilbarg2001}.

Finally, even if a PINN uses a well-posed variational formulation of a given PDE, an upper bound of the error is usually given in terms of the dual norm of the residual, which is in general not computable (see \cite{BCP} for a rigorous a posteriori error analysis for variational PINNs by a residual-type error estimator and data oscillation terms).

\subsection{Partial Differential Equations (PDEs)}
In order to explain this, let us review a situation in which the residual is related to the error a bit more in detail. For the \enquote{fruit fly of numerical analysis}, the Poisson problem on some domain $\Omega\subset\R^d$ with homogeneous Dirichlet boundary conditions, a well-posed PDE operator $\cB: H^1_0(\Omega)\to H^{-1}(\Omega)$ is induced by the bilinear form $b(u,v):=(\nabla u, \nabla v)_{L_2(\Omega)}$, where $\langle \cB u,v\rangle_{H^{-1}(\Omega) \times H^1_0(\Omega)} := b(u,v)$ and $\langle \cdot,\cdot\rangle_{H^{-1}(\Omega) \times H^1_0(\Omega)}$ denotes the duality pairing of the Sobolev space $H^1_0(\Omega)$ and its dual $H^{-1}(\Omega)$ induced by the $L_2$-inner product. The exact solution is $u=\cB^{-1}f$, where $f\in H^{-1}(\Omega)$ is a given right-hand side, i.e., $b(u,v)=\langle f,v\rangle_{H^{-1}(\Omega) \times H^1_0(\Omega)}$ for all $v\in H^1_0(\Omega)$. 
Then, coercivity implies that the error $u- u^\delta$ is bounded by the residual $r^\delta$ as follows
\begin{align*}
	\| u- u^\delta\|_{H^{1}(\Omega)}
	\le {\textstyle\frac{1}{\alpha}}\| f - \cB u^\delta\|_{H^{-1}(\Omega)}
	=: {\textstyle\frac{1}{\alpha}} \| r^\delta\|_{H^{-1}(\Omega)}, 
\end{align*}
where $u^\delta\in\X^\delta\subset\X:=H^1_0(\Omega)$ is the Galerkin solution in some (typically finite-dimensional) trial space $\X^\delta$ and $\alpha$ is the coercivity constant. Hence, the right-hand side of the latter equation might be used for the definition of the loss function and an upper error bound. This, however involves at least two issues,
\begin{itemize}
	\item the value of the coercivity constant $\alpha$ (or an lower bound) needs to be known, and
	\item the dual norm $\| r^\delta\|_{H^{-1}(\Omega)}$ of the residual needs to be computed. 
	\end{itemize}
The coercivity constant can be estimated by solving an eigenvalue problem, which might, however, be numerically costly. 
As for the second issue, recalling that 
$$
	\| r^\delta\|_{H^{-1}(\Omega)} = \sup\limits_{v\in H^1_0(\Omega)} \frac{\langle r^\delta,v\rangle_{H^{-1}(\Omega) \times H^1_0(\Omega)}}{\| v\|_{H^{1}(\Omega)}}, 
$$
one would need to determine a supremum over an infinite-dimensional space, needless to say that finite-dimensional approximations are known to violate the upper bound property in general.

\subsection{Parameterized Partial Differential Equations (PPDEs)}
This situation becomes even worse when considering Parameterized PDEs (PPDEs) of the form 
\begin{equation} \label{eq:ppde}
	\cB_{\mu} u_{\mu} = f_{\mu}
\end{equation}
where $\mu\in\cP\subset\R^p$ is some parameter.\footnote{We omit the technical assumptions for the dependency of $\cB_{\mu}$ and $f_{\mu}$ on $\mu$ and refer to the literature.}  A typical situation is that \eqref{eq:ppde} needs to be solved extremely often (multi-query) or extremely fast (realtime). In such situations, model order reduction is strongly needed. For elliptic as well as parabolic linear PPDEs allowing for an affine separation of parameters and variables, it is known that linear projection-based methods, e.g. the reduced basis method (RBM) \cite{Haasdonk:RB, Quarteroni2015, RozzaRB} are working well. On the other hand, for transport- or wave-type problems it has been proven that the Kolmogorov $N$-width decay might be poor\,\cite{OR16, greif2019decay}, such that linear model reduction techniques are bound to fail and nonlinear methods are needed. 

Since (nonlinear) model reduction can also be realized by an application of an appropriate autoencoder, it appears attractive to use PINNs also for that purpose. However, the above obstructions are even more severe here, since (a) the coercivity (or more general, the inf-sup) constant depends on the parameter $\mu$ and (b) the computation of the dual norm of the residual should be performed in realtime, i.e., extremely fast.

This is the starting point of this paper in which we suggest  to use adaptively computed wavelet expansions in order to compute the dual norm of the residual.

\subsection{Ultra-weak formulations of PDEs}
A second issue in this paper is to utilize ultra-weak formulations of PDEs which allow for an error-residual \emph{identity}, in particular without the coercivity (or more general, the inf-sup) constant, \cite{Dahmen2012, Brunken2019, henning2022ultraweak}. Such formulations are particularly appropriate for transport and wave-type problems, which we will consider in connection to PINNs in future research. Here, we review ultra-weak formulations of elliptic PDEs and show how they can be used to define certified PINNs.

Ultra-weak formulations yield an isomorphism $\cB_\mu: X:=L_2(\Omega)\to Y_\mu$\footnote{We use the symbol $\X$ for the trial and test space in the standard variational form, $X$, $Y_\mu$ for trial and test space in the ultra-weak form.} with inf-sup constant\footnote{Given in terms of the operator norm of the inverse, i.e., $\beta_\mu=\| \cB_\mu^{-1}\|^{-1}$.} being unity, so that the dual norm of the residual w.r.t.\ a suitable Petrov-Galerkin approximation $u_{\mu}^\delta\in X^\delta_\mu$ to be detailed below, namely
\begin{equation} \label{eq:residual}
	\Vert r_{\mu}^\delta \Vert_{Y_\mu'} 
		:= \Vert f_{\mu} - \cB_{\mu} u_{\mu}^\delta \Vert_{Y_\mu'}
\end{equation}
gives us not only a-posteriori information about the error of an approximation, but is also a candidate for a loss function to train the PINN.  

To evaluate the dual norm \eqref{eq:residual}, we expand the residual in a wavelet basis and exploit the norm equivalences in scales of Sobolev spaces (in particular of negative order) of these. The loss function to train the PINN is then defined as a sum of dual norms for a finite sample set $\mathcal{S}_\mu \subset \mathcal{P}$ of parameters $\mu$. Due to the fact that essential boundary and initial conditions are encoded in the residual, there is no need for a penalty term as for standard PINNs. As a result of the training process, the PINN is a nonlinear approximation of the mapping $\mu \mapsto \cB^{-1}_{\mu}f_{\mu}$, resulting in a nonlinear model reduction technique.

\subsection{Organization of the material}
The remainder of the paper is organized as follows. In Section \ref{Sec:VarPDE}, we briefly recall the crucial properties of different variational formulations in a general setting. The characterization of Sobolev spaces in terms of biorthogonal wavelets and the fast wavelet transform are described in Section \ref{Sec:Wavelets}. 
Next, we review the main ingredients of PINNs in Section \ref{Sec:PINN}. In Section \ref{sec:waveletPINN}, we derive the wavelet-based loss function to train the PINN and for an a posteriori error control. Here, for simplicity, we restrict ourselves to PPDEs with periodic boundary conditions. Numerical experiments are presented in Section \ref{sec:numExp} and Section \ref{sec:summary} is devoted to a summary and an outlook.

\section{Variational formulations for (parameterized) PDEs}\label{Sec:VarPDE}

Let $\Omega \subset \mathbb{R}^{d}$ be an open and bounded domain.\footnote{$\Omega$ might also be a domain in space and time, e.g. $\Omega=(0,T)\times\Omega_{\textrm{space}}$ for $T>0$ and $\Omega_{\textrm{space}}\subset\mathbb{R}^{\tilde{d}}$, $d=\tilde{d}+1$.} We denote the parameter by $\mu \in \mathcal{P} \subset \mathbb{R}^{p}$, the parameter set is $\mathcal{P}$. Let $\cB^{\circ}_{\mu}$ denote the \emph{classical} PPDE operator (the ${}^\circ$ indicates the classical form) defined point-wise on $\Omega$. By $\text{dom}(\cB^{\circ}_{\mu})$ we denote the \emph{classical} domain. 

\subsection{Classical form}
Essential boundary and initial conditions are encoded in the definition, so that the classical form amounts finding a function $u\in \text{dom}(\cB^{\circ}_{\mu})$ such that $\cB^{\circ}_{\mu} u_\mu(x) = f_\mu(x)$ point-wise for all $x$ and for some given right-hand side $f_\mu\in C(\overline\Omega)$.

\begin{Example}[The fruit fly]
	As a guiding example, we detail the (non-parametric) Poisson problem with periodic boundary conditions, but stress the fact that ultra-weak formulations are particularly suited for problems of transport and wave type. However, the elliptic case is particularly instructive for PINNs as it allows comparisons, so that we devote other type of problems to forthcoming papers.
	
	For $\Omega := (0,1)^{d}$, we consider the non-parametric operator $B^\circ u := -\Delta u$ along with periodic boundary conditions $u|_{x_{i}=0} = u|_{x_{i}=1}$ and $(\partial_{x_{i}} u)|_{x_{i}=0} = (\partial_{x_{i}} u)|_{x_{i}=1}$ on the boundary $\partial \Omega$, $i=1,...,d$. For a periodic right-hand side $f\in C_{\textrm{per}}(\Omega)$, there exists a unique classical solution $u\in C^1_{\textrm{per}}(\Omega)\cap C^2(\overline\Omega) \cap L_{2,0}(\Omega)$, where $L_{2,0}(\Omega):= \lbrace v \in L_{2}(\Omega): \int_{\Omega} v(x)\, dx = 0 \rbrace$.
\hfill$\triangle$
\end{Example}

\subsection{Variational formulations}

The standard variational form results from performing integration by parts and taking essential boundary (and initial) conditions into account. This is described in standard textbooks on PDEs.

\addtocounter{example}{-1}
\begin{Example}[The fruit fly, continued]
	The standard variational form is based upon the periodic Sobolev space $\X:=H^1_{\textrm{per}}(\Omega):=\mathop{clos}_{\|\cdot\|_{H^1(\Omega)}}(C^\infty_{\textrm{per}}(\Omega))$ and the bilinear form $b(u,v):=(\nabla u, \nabla v)_{L_2(\Omega)}$. The Lax-Milgram theorem ensures that for any $f\in\X'$ there exists a unique $u\in\X$ such that $b(u,v):=\langle f,v\rangle_{\X' \times \X}$ for all $v\in\X$.
	\hfill$\triangle$
\end{Example}

\subsection{Ultra-weak variational formulation}

In order to define the ultra-weak form, we consider the \emph{classical} adjoint operator $\cB^{\circ,*}_{\mu}$ of $\cB^{\circ}_{\mu}$ defined by
\begin{align*}
	(\cB^{\circ}_{\mu}u, v)_{L_{2}(\Omega)} 
	= (u, \cB^{\circ,*}_{\mu} v)_{L_{2}(\Omega)}, \quad u,v \in C_{0}^{\infty}(\Omega).
\end{align*} 
Due to \cite{Dahmen2012}, the following assumptions turn out to be crucial for the subsequent analysis:
	\begin{enumerate}
		\item[(B*1)] $\cB_\mu^{\circ,*}$ is injective on $\text{dom}(\cB_\mu^{\circ,*})$, which is assumed to be dense in $L_2(\Omega)$.
		\item[(B*2)] The range $R(\cB_\mu^{\circ,*}):= \cB_\mu^{\circ,*}[\text{dom}(\cB_\mu^{\circ,*})]$ is dense in $L_2(\Omega)$.
	\end{enumerate}
If (B*1) holds, the quantity
\begin{align*}
	\Vvert{v}\Vvert_\mu := \| \cB_\mu^{\circ,*}v\|_{L_2(\Omega)}
\end{align*}
is a norm on $\text{dom}(\cB_\mu^{\circ,*})$ and we define $Y_\mu := \text{clos}_{\Vvert{\cdot}\Vvert_\mu} \{\text{dom}(\cB_\mu^{\circ,*})\}$, which is a Hilbert space with inner product $(w,v)_\mu:= (\cB_\mu^*w, \cB_\mu^*v)_{L_2(\Omega)}$ and induced norm $\Vvert\cdot\Vvert_\mu$, where $\cB^*_\mu:Y_\mu\to L_2(\Omega)$ denotes the continuous extension of $\cB_\mu^{\circ,*}$ from $\text{dom}(\cB_\mu^{\circ,*})$ to $Y_\mu$.

\begin{Proposition}[{\cite[Prop.\ 2.1]{Dahmen2012}}]\label{theorem1}
	Let assumptions {\normalfont{(B*1)}} and {\normalfont{(B*2)}} hold and let $f_{\mu} \in (Y_{\mu})'$ be given, where $(Y_{\mu})'$ denotes the dual space of $Y_{\mu}$. Then, there is a unique $u_{\mu} \in L_{2}(\Omega)$ such that
	\begin{equation} \label{eq:ultraweakform}
		b_{\mu}(u_{\mu},v) 
		:= (u_{\mu},\cB^{*}_{\mu}v)_{L_{2}(\Omega)} 
		= f_{\mu}(v) \quad \forall v \in Y_{\mu}
	\end{equation}
	and the respective inf-sup and continuity constants are unity, i.e.,
	\begin{align}\label{eq:infsup1}
		\beta_\mu := 
		\inf_{u\in L_2(\Omega)} \sup_{v\in Y_\mu}  \frac{ b_\mu(u,v)}{\| u\|_{L_2(\Omega)}\, \Vvert{v}\Vvert_\mu} 
		= \sup_{u\in L_2(\Omega)} \sup_{v\in Y_\mu}  \frac{ b_\mu(u,v)}{\| u\|_{L_2(\Omega)}\, \Vvert{v}\Vvert_\mu} 
		= 1.
	\end{align}
\end{Proposition}

The equality in \eqref{eq:infsup1} has a consequence which is particularly relevant for this paper.  
Let $u_\mu^\delta \in L_{2}(\Omega)$ be some approximation of the exact solution $u_{\mu} \in L_{2}(\Omega)$ of \eqref{eq:ultraweakform}, e.g.\  a Petrov-Galerkin projection as suggested in \cite{Brunken2019} or a NN approximation. The optimal inf-sup stability \eqref{eq:infsup1} yields an optimal error/residual relation, namely the \emph{identity}
\begin{equation} \label{eq:errorresidual}
	\Vert e_{\mu}^\delta\Vert_{L_{2}(\Omega)} 
	:= \Vert u_{\mu}-u_\mu^\delta \Vert_{L_{2}(\Omega)} 
	= \Vert f_{\mu} - \cB_\mu u_\mu^\delta \Vert_{(Y_{\mu})'} 
	=: \Vert r_{\mu}^\delta \Vert_{(Y_{\mu})'},
\end{equation}
where $r_{\mu}^\delta$ is the residual. Here, $(Y_{\mu})'$ denotes the dual space of $Y_{\mu}$ as above.   
The error/residual identity in \eqref{eq:errorresidual} suggests to use the dual norm of the residual as a surrogate for the unknown error. However, even if we would be able to compute the residual $r_{\mu}^\delta$ (which would require the application of the \emph{exact} operator $\cB_\mu$), the computation of the dual norm is at least not straightforward since it requires to determine the supremum over the infinite-dimensional space $Y_{\mu}$. Within the RBM context, the supremum is taken only for a \enquote{truth} discretization instead of the full space. Such an approximation is appropriate for the RBM, since its computation is done via \enquote{truth} (i.e., computationally expansive) Riesz representations which can be precomputed in an offline stage. In the online stage, those precomputed quantities can highly efficiently be combined due to the assumed affine decomposition (i.e., the separation of parameter and physical variables). This, however, at least to the best of our knowledge, can not be done by using a PINN for the approximation of the solution of a PPDE. This is the reason why we use wavelet methods instead. 

\addtocounter{example}{-1}
\begin{Example}[The fruit fly, continued]
	For our guiding example above we have $X=L_{2,0}(\Omega)$ and the bilinear form reads $b(u,v):= - (u,\Delta v)_{L_2(\Omega)}$, so that $Y=H_{\textrm{per}}^{2}(\Omega)$.
	\hfill$\triangle$
\end{Example}

\section{Biorthogonal wavelets and the fast wavelet transformation}
\label{Sec:Wavelets}
Wavelets have extensively been used in signal and image analysis, but also for the numerical solution of PDEs and integral equations. In this paper, we will not use wavelet bases for numerical simulation, but (only) for determining the dual norm of the residual. The background is the fact that biorthogonal wavelets allow the characterization of scales of Sobolev spaces. Since $(Y_{\mu})'$ is usually (at least contained in) a Sobolev space with \emph{negative} index, we use a wavelet expansion of the residual combined with the norm equivalence to compute the dual norm up to any desired accuracy. In order to explain the proposed method, we will briefly review the main ingredients which we will need. For more details, we refer e.g.\ to \cite{Urban2008} and references therein. 

Even though there are wavelet bases on general domains available, we restrict ourselves to periodic wavelets on cubes in order to explain the main ingredients of our approach without introducing additional technicalities. In future research, we plan to consider extensions of the residual from some  domain $\Omega$ to a cube $\Omega\subset\Box$ and then use the method introduced here to compute the dual norm of the residual. Hence, we restrict ourselves here to PPDEs on the cube $\Omega=(0,1)^d$ with periodic boundary conditions in order to use periodic wavelets. A corresponding wavelet basis can  easily be formed by tensor products of univariate periodic wavelets on $(0,1)$. Periodic wavelets in turn arise by a simple periodization from wavelets on the real line, which we will describe next.

\subsection{Wavelets on the real line}
For any function $g: \mathbb{R} \rightarrow \mathbb{R}$, we define its \emph{scaled} and \emph{shifted} version for $j,k\in\mathbb{Z}$ by
\begin{equation}\label{eq:gjk}
	g_{[j,k]}(x) := 2^{j/2} g(2^{j} x - k), \quad x \in \mathbb{R},
\end{equation}
where $j$ is referred to as the \emph{level} and $k$ indicates the location in space (the shift). The normalization by $2^{j/2}$ is done in order to preserve the $L_2$-norm, i.e., $\| g_{[j,k]}\|_{L_2(\mathbb{R})} = \| g\|_{L_2(\mathbb{R})}$. 

Biorthogonal wavelets are typically constructed with the aid of a \emph{Multiresolution Analysis (MRA)} induced by compactly supported \emph{dual scaling functions} $\varphi$, $\tilde\varphi$. \enquote{Dual} means that  $(\varphi,\tilde\varphi(\cdot-k))_{L_2(\mathbb{R}}=\delta_{0,k}$ for all $k\in\mathbb{Z}$. A function is called \emph{scaling function}, if it is \emph{refinable}, which means here that the following two-scale relations hold
\begin{align}\label{eq:refinement}
	\varphi(x) = \sum_{k=0}^{N} a_k\, \varphi(2x-k),
	\qquad
	\tilde\varphi(x) = \sum_{k=0}^{\tilde N} \tilde a_k\, \tilde\varphi(2x-k),
	\qquad x\in\mathbb{R},
\end{align}
with real \emph{mask coefficients} $a_k$ and $\tilde a_k$ with finite $N$, $\tilde{N}\in\mathbb{N}$.\footnote{Recall that we assume that $\varphi$ and $\tilde\varphi$ are compactly supported.} Their shifted and scaled versions, i.e., $\Phi_{j}^{\mathbb{R}} := \lbrace \varphi_{[j,k]}: k \in \mathbb{Z} \rbrace$ and $\tilde{\Phi}_{j}^{\mathbb{R}} := \lbrace \tilde{\varphi}_{[j,k]}: k \in \mathbb{Z} \rbrace$ span the MRA spaces $S_j^{\mathbb{R}}=\textrm{span}(\Phi_j^{\mathbb{R}})$, $\tilde S_j^{\mathbb{R}}=\textrm{span}(\tilde\Phi_j^{\mathbb{R}})$. Then, $\lbrace S_{j}^{\mathbb{R}} \rbrace_{j \in \mathbb{Z}}$ and $\lbrace \tilde{S}_{j}^{\mathbb{R}} \rbrace_{j \in \mathbb{Z}}$ are called \emph{primal} and \emph{dual} MRAs, respectively. Examples include cardinal B-splines for $\varphi$ and dual scaling functions $\tilde\varphi$ introduced in \cite{CDF}.

Biorthogonal wavelets (on the real line) are constructed via detail spaces $W_{j}^{\mathbb{R}}$ and $\tilde{W}_{j}^{\mathbb{R}}$ defined by
\begin{align*}
	S_{j+1}^{\mathbb{R}} = S_{j}^{\mathbb{R}} \oplus W_{j}^{\mathbb{R}}, \quad W_{j}^{\mathbb{R}} \perp \tilde{S}_{j}^{\mathbb{R}}, 
	&\quad&
	\tilde{S}_{j+1}^{\mathbb{R}} = \tilde{S}_{j}^{\mathbb{R}} \oplus \tilde{W}_{j}^{\mathbb{R}}, \quad \tilde{W}_{j}^{\mathbb{R}} \perp S_{j}^{\mathbb{R}}.
\end{align*}
One then seeks (primal and dual) \emph{mother wavelets} $\psi$ and $\tilde\psi$ such that $\Psi_{j}^{\mathbb{R}} := \lbrace \psi_{[j,k]}: k \in \mathbb{Z} \rbrace$ and $\tilde{\Psi}_{j}^{\mathbb{R}} := \lbrace \tilde{\psi}_{[j,k]}: k \in \mathbb{Z} \rbrace$ are Riesz bases for the detail spaces $W_j^{\mathbb{R}}$ and $\tilde{W}_j^{\mathbb{R}}$, respectively. The collections $\Psi^{\mathbb{R}} := \bigcup_{j \in \mathbb{Z}} \Psi_{j}^{\mathbb{R}}$ and $\tilde{\Psi}^{\mathbb{R}} := \bigcup_{j \in \mathbb{Z}} \tilde{\Psi}_{j}^{\mathbb{R}}$ are called (primal and dual) wavelet bases and form biorthogonal Riesz bases for $L_{2}(\mathbb{R})$. The (compactly supported) mother wavelets can be written in terms of the scaling functions as
\begin{align}\label{eq:wavelet}
	\psi(x) = \sum_{k=0}^{M} b_k\, \varphi(2x-k),
	\qquad
	\tilde\psi(x) = \sum_{k=0}^{\tilde M} \tilde b_k\, \tilde\varphi(2x-k),
	\qquad x\in\mathbb{R},
\end{align}
and the \emph{wavelet mask coefficients} $b_k$, $\tilde b_k$ are given in terms of the mask coefficients in \eqref{eq:refinement}, \cite{CDF}. 

\subsection{Periodization}
For any scaled and shifted version $g_{[j,k]}$ in \eqref{eq:gjk}, we define its \emph{periodized version} as
\begin{equation*}
	g_{j,k} := \sum_{m\in\mathbb{Z}} g_{[j,k]}(\cdot-m)\,\chi_{[0,1]},
\end{equation*}
where $\chi_{[0,1]}$ denotes the characteristic function on the unit interval, so that the above sum over $m$ is in fact a finite one as long as $g$ has compact support. Moreover, it is readily seen that the above definition gives rise to different functions for $k=0,...,2^{j}-1$.\footnote{In fact, e.g.\ $g_{j,0}\equiv g_{j,2^j}$.}

Then, the periodic MRA spaces $S_j=\textrm{span}(\Phi_j)$, $\tilde S_j=\textrm{span}(\tilde\Phi_j)$ are induced by $\Phi_{j} := \lbrace \varphi_{j,k}: k =0,...,2^j-1 \rbrace$ and $\tilde{\Phi}_{j} := \lbrace \tilde{\varphi}_{j,k}: k=0,...,2^j-1 \rbrace$, respectively, and the biorthogonal wavelet spaces read $\Psi_{j} := \lbrace \psi_{j,k}: k =0,...,2^j-1 \rbrace$ as well as $\tilde{\Psi}_{j} := \lbrace \tilde{\psi}_{j,k}: k=0,...,2^j-1 \rbrace$, such that the detail spaces read $W_j=\textrm{span}(\Psi_j)$, $\tilde W_j=\textrm{span}(\tilde\Psi_j)$ and satisfy $S_{j+1} = S_{j} \oplus W_{j}$, $W_{j} \perp \tilde{S}_{j}$ as well as $\tilde{S}_{j+1} = \tilde{S}_{j} \oplus \tilde{W}_{j}$, $\tilde{W}_{j} \perp S_{j}$.

\subsection{The Fast Wavelet Transform (FWT)} 
Assume that we are given a function $r_J\in \tilde S_J$\footnote{In fact, we will apply this later for a sufficiently accurate approximation of the residual.}, i.e., we are given coefficients $c_{J,k}\in\R$, $k=0,...,2^J-1$, defining the function 
\begin{equation*}
	r_J = \sum_{k=0}^{2^J-1} c_{J,k}\, \tilde\varphi_{J,k}.
\end{equation*}
Recalling that $\tilde S_{j} = \tilde S_{j-1} \oplus \tilde W_{j-1}$, we can decompose $r_J$ into its pieces in $\tilde S_{J-1}$ and $\tilde W_{J-1}$, i.e.,
\begin{equation*}
	r_J  = \sum_{k=0}^{2^{J-1}-1} c_{J-1,k}\  \tilde\varphi_{J-1,k} 
		+ \sum_{k=0}^{2^{J-1}-1} d_{J-1,k}\  \tilde\psi_{J-1,k},
\end{equation*}
where the coarse level and wavelets coefficients are given by
\begin{align*}
	c_{J-1,k} 
	&= (r_{J}, {\varphi}_{J-1,k})_{L_{2}(\mathbb{R})} 
	= {\textstyle{\frac{1}{\sqrt{2}}}} \sum_{m=0}^{2^J-1} c_{J,m}\, {a}_{m-2k}
	\quad\text{ and }\\
	d_{J-1,k} 
	&
	= {\textstyle{\frac{1}{\sqrt{2}}}}  \sum_{m=0}^{2^J-1}  c_{J,m}\, {b}_{m-2k},
\end{align*}
where $a_k$, $b_k$ are the primal mask and wavelet coefficients introduced in \eqref{eq:refinement} and \eqref{eq:wavelet}, respectively. Iterating this over the levels, we obtain the \emph{Fast Wavelet Transform} (FWT)
\begin{equation*}
	\mathbf{T}_{J}: \mathbf{c}_{J} \mapsto (\mathbf{c}_{0}, \mathbf{d}_{0},...,\mathbf{d}_{J-1}),
	\qquad
	\mathbf{T}_{J}: \mathbb{R}^{2^J}\to\mathbb{R}^{2^J},
\end{equation*}
where $\mathbf{c}_{j}:=(c_{j,k})_{k=0,...,2^j-1}$ and $\mathbf{d}_{j}:=(d_{j,k})_{k=0,...,2^j-1}$. It is remarkable that the FWT and its inverse can be performed in \emph{linear} complexity, i.e., $\mathcal{O}(2^J)$ w.r.t.\ the input vector length. For simplicity, we abbreviate $\mathbf{d}_{-1}:=\mathbf{c}_{0}$, i.e., $d_{-1,k}:=c_{0,k}$.

\subsection{Characterization of Sobolev spaces.} 
Under certain (regularity and approximation) properties, it is known that biorthogonal wavelets allow for a characterization of a whole scale of Sobolev spaces, i.e.
\begin{equation} \label{eq:normequiv}
	\bigg\lVert \sum_{j\ge -1}\sum_{k=0}^{2^j-1} d_{j,k}\, \tilde\psi_{j,k} \bigg\rVert_{H^{\sigma}(\mathbb{R})}^{2} 
	\sim 
	\sum_{j\ge -1}\sum_{k=0}^{2^j-1} 2^{2\sigma j} \vert d_{j,k} \vert^{2}, 
	\quad \forall \sigma \in (-{\gamma}, \tilde\gamma),
\end{equation}
where $\gamma, \tilde{\gamma}>0$ are constants depending on regularity and compression properties of $\varphi$ and $\tilde\varphi$. Here $A \sim B$ abbreviates the existence of constants $0 < c \le C < \infty$ such that $cA \le B \le CA$. The specific constants in \eqref{eq:normequiv} can be estimated by solving certain eigenvalue problems, which we also did in our subsequent numerical experiments.

\subsection{Dual norm of the residual}\label{Sec:DualNormResidual} 
Note, that \eqref{eq:normequiv} involves also Sobolev spaces of \emph{negative} order, which we use for determining the dual norm of the residual, $\Vert r_{\mu}^\delta \Vert_{Y'}$ or $\Vert r_{\mu}^\delta \Vert_{(Y_{\mu})'}$. We shall detail this norm equivalence for the variational formulations under consideration.

\subsubsection{Standard variational form}
In the standard variational formulation, the residual is in a standard periodic Sobolev space of negative order, i.e., $\sigma$ in \eqref{eq:normequiv} is the negative of half of the order of the PDE.

\addtocounter{example}{-1}
\begin{Example}[The fruit fly, continued]
	We need $\X'$ for $\X:=H^1_{\textrm{per}}(\Omega)$, i.e., $\sigma=-1$.
	\hfill$\triangle$
\end{Example}

\subsubsection{Ultra-weak form}
The ultra-weak formulation arises from putting \emph{all} derivatives involved in the definition of PPDE onto the test functions by means of integration by parts. Hence, $Y_{\mu}$ is a problem-dependent smoothness space. In all examples, that we have considered, the dual space can be associated to a Sobolev space with negative index, i.e., 
\begin{equation*}
(Y_{\mu})' \cong (H^{\sigma}_{\textrm{per}}(\Omega))',
\qquad
H^{\sigma}_{\textrm{per}}(\Omega) := \text{clos}_{\Vert \cdot \Vert_{H^{\sigma}}}(C^{\infty}_{\textrm{per}}(\Omega)),
\end{equation*}
with $Y_{\mu}$ containing periodic test functions. The exponent $\sigma$ is given by the order of the PPDE. 

\addtocounter{example}{-1}
\begin{Example}[The fruit fly, continued]
	In the ultra-weak form we get $(Y_{\mu})' \cong (H^{2}_{\textrm{per}}(\Omega))'$. At least for $d=1$, the norm $\Vvert{\cdot}\Vvert_\mu$ can be shown to be equivalent to the standard $H^{2}(\Omega)$-norm (this means that we get $\sigma=-2$) with a $\mu$-dependent constant, see \S\ref{Sec:Num2UWF} below.
	\hfill$\triangle$
\end{Example}

From now on, for simplicity of the exposition, we restrict ourselves to $\Omega=(0,1)$.

\subsubsection{Wavelet coefficients}
For using the norm equivalence in \eqref{eq:normequiv} provided by appropriate biorthogonal wavelets, we obviously need the \emph{wavelet coefficients} of the residual, namely
\begin{align}\label{eq:djk}
	d_{j,k}(r_{\mu}^\delta) := (r_{\mu}^\delta, \psi_{j,k})_{L_2(0,1)},
\end{align}
i.e., inner products of the residual with the primal wavelets.\footnote{We will later use the residual $r_\mu^\theta$ of the NN instead of $r_{\mu}^\delta$, which is the reason why we make the dependency of the expansion coefficients on the underlying functions explicit.} We are going to show two approaches for computing the desired wavelet coefficients: one using the FWT, the other one in terms of an adaptive computation (where both approaches can also be combined).

\subsubsection{Using the FWT}
Since the primal wavelets by \eqref{eq:wavelet} are a linear combination of the primal scaling functions $\varphi(2\cdot-k)$ and $\varphi$ is a B-spline, we can in fact compute $d_{j,k}$ by using the FWT based upon the \emph{single scale coefficients} $c_{J,k}$, where
\begin{align}\label{eq:rhoJ}
	r_{\mu}^{\delta,J} 
		&:= \sum_{k=0}^{2^J-1} c_{J,k}(r_{\mu}^\delta)\, \tilde\varphi_{J,k} \in \tilde{S}_J,
	\qquad 
	c_{J,k}(r_{\mu}^\delta) = (r_{\mu}^\delta, \varphi_{J,k})_{L_2(0,1)}.
\end{align}
Of course, $r_{\mu}^{\delta,J}$ in \eqref{eq:rhoJ} is only an approximation to $r_{\mu}^\delta$. However, by increasing the level $J$, we can achieve \emph{any} desired accuracy for the error $\| r_{\mu}^\delta - r_{\mu}^{\delta,J}\|_{L_2(0,1)}$.  In fact, a standard Jackson-type inequality yields $\| r_{\mu}^\delta - r_{\mu}^{\delta,J}\|_{L_2(0,1)} \le C\, 2^{-Js} \| r_{\mu}^\delta\|_{H^s(0,1)}$. 
The corresponding wavelet coefficients are obtained by performing the FWT for the expansion in terms of $\tilde\Psi$.

Since $\varphi$ is a cardinal B-spline, the computation of the coefficients $c_{J,k}(r_{\mu}^\delta)$ is in fact straightforward by numerical quadrature as long as we can sample the residual $r_{\mu}^\delta$ point-wise. This is possible as long as we have access to the right-hand side $f_\mu$ (since the PINN provides us with point values of the operator $\cB_\mu$ applied to the PINN approximation). 

\subsubsection{Adaptive computation}
However, recalling that the residual is a functional (an element of a dual space), point-wise evaluation might not make sense.\footnote{We are particularly interested in the low regularity case, where the solution has the minimal regularity, so that the residual cannot be interpreted as a function, but must be treated as a functional.} Hence, we consider the \emph{expansion} of the residual in terms of the dual wavelet basis $\tilde\Psi$, where the primal wavelets $\Psi$ are at least in $H^{\sigma}_{\textrm{per}}(0,1)$. This can easily be achieved since the primal scaling function $\varphi$ can be chosen as a cardinal B-spline of corresponding order, so that $\psi$ is a linear combination of B-splines. Then,
\begin{align*}
	{r}^\delta_{\mu} = \sum_{j=-1}^\infty \sum_{k=0}^{2^j-1} d_{j,k}(r^\delta_{\mu})\,  \tilde\psi_{j,k},
	\quad\text{where}\quad
	d_{j,k}(r^\delta_{\mu}) = (r^\delta_{\mu}, \psi_{j,k})_{L_2(0,1)},
\end{align*}
i.e., we need the inner products of the residual with B-splines as described above. Now, we can use the norm equivalence in \eqref{eq:normequiv} to obtain
\begin{equation}\label{eq:weightedsum}
	\Vert {r}_{\mu}^\delta \Vert_{(H^{\sigma}_{\textrm{per}}(0,1))'}^2
	\sim \sum_{j=-1}^\infty \sum_{k=0}^{2^j-1} 2^{-2\sigma j} |d_{j,k}(r^\delta_{\mu})|^2.
\end{equation}
As such it is impossible to compute the right-hand side of the latter equation exactly as it involves a sum over an infinite index range. However, we can now use well-known compression properties of wavelets to approximate the infinite sum up to any desired accuracy. In fact, the following Whitney-type estimate is well-known to hold
\begin{equation}\label{eq:whitneyesti}
	\vert d_{j,k}(r^\delta_{\mu}) \vert 
	\le C\, 2^{-\sigma j} \Vert {r}_{\mu}^\delta \Vert_{\sigma,\,  \text{supp } \psi_{\lambda}}, \quad \text{for } 0 \le \sigma < \tilde{d},
\end{equation}
where $\tilde{d}$ is the number of vanishing moments\footnote{$\tilde{d}$ coincides with the \emph{order} of the dual scaling function $\tilde\varphi$ and can be chosen sufficiently large, \cite{CDF,Urban2008}.} of the wavelets. This means that wavelet coefficients are small in regions where the residual is smooth. This allows for an adaptive computation of the dual norm similar to  \cite{AliSteihUrban,AliUrban} within the framework of the RBM.

No matter how we compute the wavelet coefficients, once they are determined, we evaluate the scaled sum on the right-hand side of \eqref{eq:weightedsum} as a surrogate for the dual norm of the residual.

\section{Physics informed neural networks (PINNs) for PPDEs}
\label{Sec:PINN}
Since we aim at investigating the use of PINNs for solving PPDEs, we  briefly review the required facts on NNs and PINNs. The subsequent paragraph is based upon  \cite{Berner2021,Gribonval2021,Petersen2018}. As usual, we define a NN as a function\footnote{With a slight abuse of notation we use the symbol $\Phi$ both for the scaling function bases and for the NN, but the notation should be clear from its use.} $\Phi_a(\cdot;\theta):\mathbb{R}^{N_{0}} \rightarrow \mathbb{R}^{N_{L}}$, which is determined by the fixed architecture $a$ of the NN and the parameters $\theta$, which are subject to a training phase. The architecture $a=(N,\rho)$ depends on the dimension $N=(N_0,...,N_L)\in\mathbb{N}^{L+1}$, $L \in \mathbb{N}$ denoting the number of \emph{layers} ($N_0$ being the input and $N_L$ the output dimension) and the \emph{activation function} $\rho: \mathbb{R} \rightarrow \mathbb{R}$. Moreover, $N_l$ is the number of \emph{neurons} in layer $l=0,...,L$ and the layers for $l=1,...,L-1$ are called \emph{hidden}.

The parameters of the NN read $\theta = (W^{(l)},b^{(l)})_{l=1,...,L}$, where $W^{(l)} \in \mathbb{R}^{N_{l} \times N_{l-1}}$ are the \emph{weight matrices} and $b^{(l)} \in \mathbb{R}^{N_{l}}$, are called \emph{bias vectors}. 
The output $\Phi_a(z;\theta)$ of the NN for an input $z\in\mathbb{R}^{N_0}$ is then defined as $\Phi_{a}(z;\theta) := \Phi^{(L)}(z;\theta)$, where
	\begin{align*}
		\Phi^{(1)}(z;\theta) &= W^{(1)} z + b^{(1)}, \\
		\hat{\Phi}^{(l)}(z;\theta) &= \rho(\Phi^{(l)}(z;\theta)), \quad l=1,...,L-1, \quad \text{and} \\
		\Phi^{(l+1)}(z;\theta) &= W^{(l+1)} \hat{\Phi}^{(l)}(z;\theta) + b^{(l+1)}, \quad l = 1,...,L-1,
	\end{align*}
and $\rho$ is applied component-wise.

\emph{Physics informed neural networks (PINNs)} are NNs, whose loss function is particularly suited for the problem to be solved, e.g.\ it uses the residual of a PDE. We are now going to describe the main ingredients of PINNs for the specific case of a PPDE, see also  \cite{Raissi2019}. 

\subsection{Standard PINNs for the classical form of PPDEs} 
As above, let $\cB_\mu^\circ$ denote the classical PPDE operator, where we again assume that essential boundary conditions, namely $u_\mu(x)=g_\mu(x)$, $x\in\Gamma:=\partial\Omega$ for some given $g_\mu\in C(\Gamma)$, are incorporated in the definition of the operator. Then, we define the point-wise classical inner and boundary residual, respectively, for some  $v\in\textrm{dom}(\cB^\circ_\mu)$ as
\begin{align*}
	r_{\Omega}(v; x,\mu) := f_{\mu}(x) - (\cB^{\circ}_{\mu} v)(x), \,\, x\in\Omega, 
	&\qquad
	r_{\Gamma}(v; x,\mu) := g_{\mu}(x) - v(x),\,\, x\in\Gamma.
\end{align*}
We fix the architecture $a=(N,\rho)$ of the NN and construct a PINN $\Phi_a(z;\theta_{\text{PINN}})$ for the input $z=(x,\mu)\in \overline\Omega\times\mathcal{P}\subset\mathbb{R}^d\times\mathbb{R}^p$, which means that the input dimension is $N_0=d+p$. The output is the value for an approximation to $u_\mu$ at $x\in\overline\Omega$, i.e., the output dimension is thus $N_L=1$. Of course, one could also consider functions of the state, namely output functions. In that case, a primal-dual approach seems appropriate.

The loss function for a (classical) PINN\footnote{Classical in the sense that the PINN corresponds to the classical form of the PPDE.} is defined by the point-wise residuals, i.e.,
\begin{align}\label{eq:stronglosspart}
	\mathcal{R}_{\textrm{classical}} (v; x,\mu) :=
	\begin{cases}
		r_{\Omega}(v; x,\mu), 	& \text{ if } x\in\Omega, \\
		\omega_b\, r_{\Gamma}(v; x, \mu), 	& \text{ if } x\in\Gamma,
	\end{cases}
\end{align}
where $\omega_b\in\mathbb{R}^+$ is a weight factor, \cite{Wang2021}. This means that training data, e.g.\ derived from numerical simulations, are not needed since the residual should vanish for the exact classical solution. Hence, one chooses finite training sets $\mathcal{S}_{\Omega} = \lbrace x_{i} \in \Omega: i=1,...,n_{\Omega} \rbrace$, $\mathcal{S}_{\Gamma} = \lbrace x_{i} \in \Gamma: i=1,...,n_{\Gamma} \rbrace$, $\mathcal{S}_x:= \mathcal{S}_{\Omega}\cup\mathcal{S}_{\Gamma}$ and $\mathcal{S}_{\mu} = \lbrace \mu_{i} \in \mathcal{P}: i=1,...,n_{\mathcal{P}} \rbrace$. The loss function is then defined by the squared residuals
\begin{equation}\label{eq:strongloss}
	\mathcal{L}_{\textrm{classical}} (v) 
	:= \sum_{(x,\mu) \in \mathcal{S}_{x}\times\mathcal{S}_{\mu} } 
		\mathcal{R}_{\textrm{classical}} (v; x,\mu)^2.
\end{equation}
The desired parameters $\theta_{\text{PINN}}$ are determined in the training phase as
\begin{equation} \label{eq:opProblem}
	\theta_{\text{PINN}} = \arg\min_{\theta} \mathcal{L}_{\textrm{classical}} (\Phi_{a}(\cdot;\theta)).
\end{equation}
In order to compute $r_{\Omega}(\Phi_{a}(\cdot;\theta); x,\mu)$, the values of $\cB^{\circ}_{\mu}(\Phi_{a}(\cdot;\theta))$ are needed for the samples in $\mathcal{S}_{x}\times\mathcal{S}_{\mu}$. These values can be determined by backpropagation, which requires a certain regularity of the activation function, e.g.\ the rectified linear unit (ReLU) function $\rho_{\text{relu}}(x) := \max \lbrace 0, x \rbrace$ can not be used for PDE operators of order more than one.

\subsection{Standard PINNs for the (standard) variational form of PPDEs}
It is an obvious drawback of the above PINN-approach that it relies on the classical form of a PPDE, which might not be well-posed. Hence, there are attempts to define a corresponding loss function for weak-form residuals, see e.g.\,\cite{Khodayi-Mehr2020,Kharazmi2021,Kharazmi2019}. The main idea is to modify the residual terms in $\Omega$ to a variational form. To this end, let $\mathcal{S}_\Omega\subset L_2(\Omega)$ be a finite sample set (of functions). The residual of the variational form is tested by such functions, \cite{Kharazmi2019,Khodayi-Mehr2020}. The arising form\footnote{This might be an inner product or a duality pairing.} is denoted by $\langle r_\Omega(v; \cdot,\mu), \phi\rangle$ for $\phi\in\mathcal{S}_\Omega$ and \eqref{eq:strongloss} is replaced by the quantity
\begin{equation}\label{eq:weakloss}
	\mathcal{L}_{\textrm{weak}} (v) 
	:= \sum_{\mu \in \mathcal{S}_{\mu} }
		\left[\sum_{\phi \in \mathcal{S}_\Omega} 
		\langle r_\Omega(v; \cdot, \mu), \phi\rangle^2
			+ \omega_{b} \sum_{x \in \mathcal{S}_{\Gamma}} (r_{\Gamma}(v; x,\mu))^{2} 
			\right],
\end{equation}
which can then be used as a loss function replacing $\mathcal{L}_{\textrm{classical}}$ in the training \eqref{eq:opProblem}. Also the boundary residual can be enforced in a weak manner, \cite{Kharazmi2021}. However, it is known that the training does not necessarily converge, see \cite{Wang2021,Wang2022} for cases of failures and \cite{Shin2020v1, Shin2020v2,Jiao2022} for examples showing convergence. Also, the choice of an appropriate weight factor $\omega_b$ might be delicate. 

In any case, at least to the very best of our knowledge, the above described loss functions do not yield error control in the sense that the accuracy of a PINN output can be rigorously bounded.\footnote{An exception is \cite{BCP}, where an a posteriori error analysis is presented for a specific Petrov-Galerkin realization of a variational PINN.} The reasons are 
\begin{itemize}
	\item the well-posedness of the PPDE requires inf-sup stability, which needs to be ensured;
	\item in case of well-posedness, the error-residual relation involves the inverse of the inf-sup constant, \cite{XuZikatanov}. This, however, might not be accessible for a \emph{parameterized} PDE, or one would need to resort to lower bounds, \cite{SCM};
	\item in addition, the dual norm of the residual needs to be determined in realtime, see \eqref{eq:errorresidual}. 
\end{itemize}
As a consequence, both loss functions mentioned above do not yield a reliable upper bound for the error. This is the reason why we also investigate the ultra-weak form allowing an error-residual \emph{identity} (even for transport and wave-type problems) in combination with wavelets for the computation/approximation of the dual norm of the residual.

\subsection{Standard PINNs for the ultra-weak variational form of PPDEs}
We fix an architecture $a=(N,\rho)$ of a NN such that $\rho \in C(\mathbb{R})$, e.g.\ the ReLU function\footnote{Which is indeed possible for the ultra-weak form as opposed to previous form of PINNs.}. For that setting, NN-approximation results for $L_2$-functions are known, see e.g.\ \cite{Guehring2020}. Recalling that the ultra-weak solution is searched in $L_2$, our aim is to control the error (the difference between the exact ultra-weak solution $u_\mu$ in \eqref{eq:ultraweakform} and the output $\Phi_a(\cdot,\mu;\theta)$ of the NN) in the $L_2$-norm, i.e., 
\begin{align}\label{eq:optimaltarget}
	\| u_\mu -\Phi_a(\cdot,\mu;\theta)\|_{L_2(\Omega)}^2
	= \Vert f_{\mu} - \cB_\mu \Phi_a(\cdot,\mu;\theta) \Vert_{(Y_{\mu})'}^2 
	=: \Vert r_{\mu}^\theta \Vert_{(Y_{\mu})'}^2
	=: \mathcal{L}_\mu (\theta) 
\end{align}
for any input parameter $\mu$. In this sense, we call the PINN \emph{certified}. Choosing again a finite parameter sample set $\mathcal{S}_{\mu}\subset\mathcal{P}$ of parameters, the loss function takes the form
\begin{equation}\label{eq:uwloss}
	\mathcal{L}_{\textrm{ultra-weak}} (\theta) 
	:= \sum_{\mu \in \mathcal{S}_{\mu} } 
		\mathcal{L}_\mu (\theta).
\end{equation}
We propose to compute $\mathcal{L}_\mu (\theta)$ for a given parameter $\mu$ by using wavelet methods.

\section{Wavelet-based PINNs}\label{sec:waveletPINN}
We are now going to describe in detail how we use wavelet methods for the computation of a loss function for a PINN for PPDEs. To this end, we use $\Vert r_{\mu}^\theta \Vert_{\X'}$ in case of the standard variational form and $\Vert r_{\mu}^\theta \Vert_{(Y_{\mu})'}$ for the ultra-weak case within the definition of the loss function. In other words, we replace $r_{\mu}^\delta$ in \S\ref{Sec:DualNormResidual} by the NN-residual $r_{\mu}^\theta$. 

Moreover, these dual norms are determined / approximated as described above by appropriately weighted wavelet coefficients. This means that we approximate the loss function at an a priori chosen accuracy and then minimize this approximation during the training phase. We can use either the FWT or the adaptive procedure for computing the wavelet coefficients. The FWT-variant is detailed in  Algorithm\,\ref{algo:evalAlgo} (by $\|\cdot\|$ we denote the Euclidean vector norm). By choosing the regularity $\sigma$ appropriately, Algorithm\,\ref{algo:evalAlgo} can realize the standard variational version (denoted then by $\mathcal{L}_{\textrm{weak},J}$) or the ultra-weak one (yielding $\mathcal{L}_{\textrm{ultra-weak},J}$).

{\centering
\begin{minipage}{.99\linewidth}
\begin{algorithm}[H] \label{algo:evalAlgo}
	\SetAlgoLined
	\SetKwInOut{Input}{Input}\SetKwInOut{Output}{Output}
	\Input{Finite parameter sample $\mathcal{S}_{\mu} \subset \mathcal{P}$, 
		NN-residual ${r}_{\mu}^\theta$, regularity $\sigma>0$, \\
		scaling function $\varphi$, finest level $J$, wavelet transform $\mathbf{T}_{J}$}
	\Output{Training loss approximation $\mathcal{L}_{\textrm{wavelet},J} (\theta)$ of \eqref{eq:uwloss} }
	$\mathcal{L}_{\textrm{wavelet},J} (\theta) \leftarrow 0$\;
	\For{$\mu \in \mathcal{S}_{\mu} $}{
		$\mathbf{c}_{J} \leftarrow (({r}_{\mu}^\theta, \varphi_{J,k})_{L_2(\Omega)})_{k=0,...,2^J-1}$
		\hfill\tcp{single-scale coefficients}
		\label{algo:coeffsC}
		$(\mathbf{c}_{0}, \mathbf{d}_{0}, \ldots,\mathbf{d}_{J-1}) \leftarrow \mathbf{T}_{J} \mathbf{c}_{J}$  
		\hfill\tcp{fast-wavelet transform}
		\label{algo:fwt}
		$\mathcal{L}_{\textrm{wavelet},J} (\theta)\, += 
			\Vert \mathbf{c}_{0} \Vert^{2}
			+\sum\limits_{\ell = 0}^{J-1}  2^{-2\sigma\ell} \Vert \mathbf{d}_{\ell} \Vert^{2}$\; 
		\label{algo:frobnorm}
	}
	\caption{Loss function of Wavelet-FWT-based PINN}
\end{algorithm}
\end{minipage}
\par
}

Some remarks are in order:
\begin{itemize}
	\item For $\Omega=(0,1)^d$, the computational complexity is of the order $2^{Jd}$ as this is the size of the scaling function basis and the FWT is of linear complexity. Of course one could resort to known approaches in order to deal with the curse of dimensionality w.r.t.\ $d$. However, this issue is independent of our approach.
	\item The computation of $({r}_{\mu}^\theta, \varphi_{J,k})_{L_2(\Omega)}$ is done by numerical integration. This amounts point-wise evaluation of the residual. At this point, however, recall that the ultra-weak form of the PPDE yields an $L_2$-approximation. Hence, $\cB_\mu(\Phi_a(\cdot,\mu;\theta))$ may not be defined point-wise. Hence, we use integration by parts again (which is possible if $\varphi_{J,k} \in Y_{\mu}$) and compute the scaling function coefficients by
	\begin{align*}
		c_{J,k}({r}_{\mu}^\theta) 
		&
		= (f_\mu - \cB_\mu \Phi_a(\cdot,\mu;\theta), \varphi_{J,k})_{L_2(\Omega)}\\
		&= (f_\mu, \varphi_{J,k})_{L_2(\Omega)} 
			- (\Phi_a(\cdot,\mu;\theta), \cB_\mu^*\varphi_{J,k})_{L_2(\Omega)}.
	\end{align*}
	Choosing the scaling function $\varphi$ sufficiently smooth allows us to compute both terms by quadrature. For the variational variant, we use standard integration routines.
\end{itemize}

\section{Numerical experiments} \label{sec:numExp}
In this section, we report some results of our numerical experiments, which have been performed using PyTorch, \cite{NEURIPS2019_9015}. The optimization for the training phase has been done by the truncated Newton method (TNC) of the Scipy package\,\cite{2020SciPy-NMeth} using the full gradient. Since the output of the training depends on the initialization of the weights, each optimization run has been repeated five times; the lines in the plots below show the mean. Further, the termination criterion of the training phase is a fixed number of training steps.

\subsection{Poisson problem, nonparametric}
We start by a simple example for validation purposes. For $\Omega := (0,1)$, we consider the non-parametric ($\cB_\mu^\circ u \equiv \cB^\circ u$, $f_\mu \equiv f$) Poisson problem with smooth data  and periodic boundary conditions, i.e., 
\begin{equation}\label{eq:nonParamProblem}
		\cB^\circ u := - u'', \quad
		f(x) := (4 \pi^{2})\, \cos(2 \pi x), 
		\quad
		u(0) = u(1), u'(0) = u'(1).
\end{equation}
Hence, the solution is smooth. Moreover, since the analytical solution $u_{\textrm{ana-sol}}$ is known, we can use a standard NN to approximate $u_{\textrm{ana-sol}}$ by using the point-wise error for the loss function. To be precise, we consider the MSE-training\footnote{MSE: mean-squared error} with the loss function 
\begin{equation} \label{eq:mseLoss}
	\mathcal{L}_{\textrm{NN}}(\theta) 
	:= \frac{1}{\vert \mathcal{S}_{x} \vert} \sum_{x \in \mathcal{S}_{x}} 
			\vert u_{\textrm{ana-sol}}(x) - \Phi_{a}(x; \theta) \vert^{2}.
\end{equation}
The training set $\mathcal{S}_{x}$ is chosen to be the set of quadrature points used in Algorithm\,\ref{algo:evalAlgo}, so that we can perform a meaningful comparison with our wavelet-based PINN for the variational formulations.\footnote{With three quadrature points in 1D the size of the training set is $\vert \mathcal{S}_{x} \vert = 3 \cdot 2^{J}$.} 
Of course, the data of this example is particularly in favor for the classical NN approach. We do not consider a PINN using the $L_2$-norm of the residual since the NN is \enquote{best one can do with a NN}.

We shall compare this NN with the classical PINN and our wavelet-based PINN for the standard variational form as well as for the ultra-weak form. Of course, for such an elliptic problem, the ultra-weak form is somewhat non-standard, but it is instructive to investigate how the different variants of our method compare in this simple case.  We recall from the \enquote{fruit fly} example that $Y_\mu \equiv Y:= H_{\textrm{per}}^{2}(\Omega)$, $b_\mu(u,v)\equiv b(u,v):=-(u,\Delta v)_{L_2(\Omega)}$ and $X=L_{2,0}(\Omega)$. We use the $L_{2}$-orthogonal projector $Q: L_{2}(\Omega) \rightarrow L_{2,0}(\Omega)$ as a normalization. 

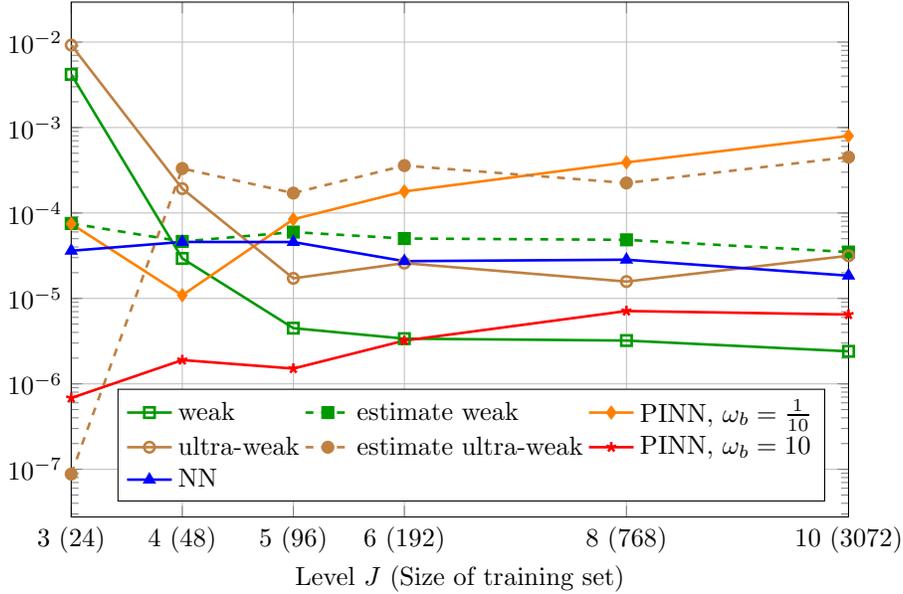
\begin{figure}[htb]
	\centering
	\begin{tikzpicture}[scale = 0.99]
	\begin{semilogyaxis}[
	legend cell align={left}, 
	xlabel=Level $J$ $\left( \text{Size of training set} \right)$,
	grid=major,
	xticklabels={$3 \ (24)$, $4 \ (48)$, $5 \ (96)$, $6 \ (192)$, $8 \ (768)$, $10 \ (3072)$},
	width=0.95\textwidth,
	height=8.5cm,
	legend pos=south east, 
	legend columns=3,
	cycle list name=black white,
	xmin=3,xmax=10,
	xtick=data]
	\addplot[mark=square,green!60!black,line width=1pt] table [x=level, y=wError] {figures/error_tc6_uw_w_pw_cl.dat};
	\addlegendentry{weak}
	\addplot[mark=square*,green!60!black,line width=1pt,dashed,mark options={solid}] table [x=level, y=wEstError] {figures/error_tc6_uw_w_pw_cl.dat};
	\addlegendentry{estimate weak}
	\addplot[mark=diamond*,orange,line width=1pt] table [x=level, y=PINNError01] {figures/error_tc6_uw_w_pw_cl.dat};
	\addlegendentry{PINN, $\omega_{b}=\frac{1}{10}$}
	\addplot[mark=o,brown,line width=1pt] table [x=level, y=uwError] {figures/error_tc6_uw_w_pw_cl.dat};
	\addlegendentry{ultra-weak}
	\addplot[mark=*,brown,line width=1pt,dashed,mark options={solid}] table [x=level, y=EstError] {figures/error_tc6_uw_w_pw_cl.dat};
	\addlegendentry{estimate ultra-weak}
	\addplot[mark=star,red,line width=1pt] table [x=level, y=PINNError10] {figures/error_tc6_uw_w_pw_cl.dat};
	\addlegendentry{PINN, $\omega_{b}=10$}
	\addplot[mark=triangle*,blue,line width=1pt] table [x=level, y=MSEError] {figures/error_tc6_uw_w_pw_cl.dat};
	\addlegendentry{NN}
	\end{semilogyaxis}
	\end{tikzpicture} 
	\caption{\label{fig:l2comparison}Comparison of the $L_{2}$-error after training with different loss functions for\,\eqref{eq:nonParamProblem}. All NNs have the same architecture $a = (N,\rho)$, where $N = (1, 64, 64, 64, 1)$ and $\rho(x) = \tanh(x)$. The level $J$ refers to the finest level in Algorithm\,\ref{algo:evalAlgo} (in brackets the resulting size of the training set).}
\end{figure}

The results are depicted in Figure\,\ref{fig:l2comparison}. We show different quantities over the level $J$ used for the quadrature points and in the FWT-computation of the dual norm of the residual. All norms are measured in $L_2(\Omega)$ in order to get a fair comparison (even though the $H^1$-norm would be more appropriate for the standard weak form, see the remarks and also \S\ref{Sec:NumPPDE} below). Since the training depends on the initialization, we display the mean. Here, the analytical solution is known, hence we can compute exact errors, which we did for different choices of the (PI)NN:
\begin{itemize}
	\item green: wavelet-based PINN using the standard weak form / error estimate;
	\item brown: wavelet-based PINN using the ultra-weak form / error estimate;
	\item orange/red: classical PINN \eqref{eq:stronglosspart} for two different choices of the weight factor $\omega_b$ for the boundary conditions;
	\item blue: NN for approximating the analytical solution with MSE-training. 
\end{itemize}

We see that most curves stagnate as $J$ increases (we comment below on those cases, where this does not seem to be the case). This means that already small levels $J$ are sufficient for the computation of the wavelet-based loss function.\footnote{This was expected since the solution and the right-hand are smooth so that the residual is also smooth implying a fast convergence of the single-scale approximation as $J\to\infty$.} 
Moreover, we observe that the classical PINN with $\omega_b=0.1$ (orange) in \eqref{eq:stronglosspart} is too weak in the sense that the resulting error is large even though the training error is small. Moreover, the error does not seem to stagnate for $J=10$, which is due to the poor choice of $\omega_b$. On the other hand, for $\omega_b=10$ (red), we obtain good results. We stress the fact that $\omega_b$ depends on the size of $\mathcal{S}_{\Omega}$ as well as $\mathcal{S}_{\Gamma}$ and a \enquote{good} choice is problem-dependent, so that typically it is not known in practice. This shows that classical PINNs might not be certified even for such a simple example.

Next, we observe that the wavelet-based error bounds (dashed lines) are in fact strict upper bounds of the error already for $J\ge 4$ for both variational forms. However, they are by one order of magnitude larger than the exact error (solid lines with the same color). The reason is that we did not optimize the constant in the norm equivalence  \eqref{eq:weightedsum}. Since we see that the slopes of the error estimate and the exact error are basically parallel, optimizing this constant will lead to a significant improvement of the effectivity, see also \S\ref{Sec:NumPPDE} below.

The weak variational form (green) shows the best results, which is not completely surprising. This formulation uses the $H^{-1}$ norm of the residual, i.e., the solution in $H^1$, which means that the gradient is also controlled. Hence, the error is measured in a stronger norm, which results in the observation that the $L_2$-norm is smaller than for the NN and the ultra-weak case. 

We also compare the wavelet-based ultra-weak PINN (brown) with the classical NN-approximation of the exact solution (blue). For $J\ge 4$, both curves are quite close. Hence, the ultra-weak form (even though using low regularity) provides results that are comparable to the NN-approximation of the exact solution, and allows for error certification.

Finally, we comment on the observation that the curves for the ultra-weak case (brown) do not seem to stagnate for $J=10$. The reason is that we fixed the number of iterations in the training phase for comparison purposes. It is known (see e.g.\ \cite{henning2022ultraweak}) that the condition numbers in the ultra-weak case are significantly larger than for the standard variational form. Hence, we would need to adjust the number of iterations in the training to see the stagnation w.r.t.\ the level.

\textbf{CPU times.} In Table \ref{tab1}, we monitor the CPU times for the computation of the loss function in the weak, ultra-weak, classical and NN-approximation case. Of course, if the formula for the exact solution is known, the evaluation of $\mathcal{L}_{\textrm{NN}}$ can be performed highly efficient. However, also the evaluation of ultra-weak loss function is quite fast, in particular when comparing it to the classical-form PINN. The CPU time for the standard weak form is halfway between ultra-weak and classical form.
\begin{table}[htb] 
	\centering
	\small
	\caption{CPU times [seconds] for the evaluation of the loss function for different levels $J$ and number of sample points, respectively.\label{tab1}}
	\begin{tabular}{rcccc}
		\toprule
		\textbf{Level} $J$ & $\mathcal{L}_{\textrm{ultra-weak},J}$  [s] & $\mathcal{L}_{\textrm{weak},J}$ [s] & $\mathcal{L}_{\textrm{classical}}$ [s] &  $\mathcal{L}_{\textrm{NN}}$ [s] \\
		\midrule
			3 & 0.00467 & 0.00916 & 0.00754 & 0.00111 \\
			4 & 0.00638 & 0.01019 & 0.00736 & 0.00110 \\
			5 & 0.00800 & 0.01172 & 0.00757 & 0.00114 \\
			6 & 0.00953 & 0.01394 & 0.00764 & 0.00140 \\
			8 & 0.01249 & 0.01749 & 0.00976 & 0.00125 \\
			10 & 0.01640 & 0.02361 & 0.01777 & 0.00189 \\
			12 & 0.02623 & 0.04393 & 0.05687 & 0.00573 \\
			14 & 0.05519 & 0.11520 & 0.21085 & 0.02575 \\
		\bottomrule
	\end{tabular}
\end{table}

\subsection{Parameterized diffusion-reaction problem}\label{Sec:NumPPDE}
Next, we consider a PPDE on $\Omega=(0,1)$ where the variable diffusion coefficient acts as parameter. The right-hand side is chosen as in \eqref{eq:nonParamProblem}, i.e., independent of the parameter. The remaining data reads
\begin{equation}\label{eq:ParamProblem}
		\cB_\mu^\circ u(x) := - \left( A(x;\mu) u'(x;\mu) \right)' + u(x;\mu), 
		\qquad
		u(0) = u(1),\, u'(0) = u'(1).
\end{equation}
where the parameter-dependent coefficient function $A(\mu)=A(\cdot;\mu): \Omega\to\R$ is defined as a continuous piecewise affine function ($p$ is the dimension of the parameter) for $d=1$ by
\begin{equation*}
	A(x;\mu) := \sum\limits_{i = 1}^{p+1} 
		\left\{\frac{\mu_{i}-\mu_{i-1}}{\xi_{i}-\xi_{i-1}} x 
		+ \frac{\mu_{i-1}\xi_{i} - \mu_{i} \xi_{i-1}}{\xi_{i}-\xi_{i-1}}\right\} \chi_{[\xi_{i-1},\xi_i)}(x), \quad \mu_{0} = \mu_{p+1} = 1,
\end{equation*}
where $\chi_{[a,b)}$ denotes the characteristic function of the interval $[a,b)$ and $\xi_i:=\frac{i}{p+1}$ are the locations where $A(\xi_{i};\mu) = \mu_{i}$ holds. For our experiments, we used $p=2$ and $\cP := [\frac18, 2]^{2} \subset \mathbb{R}^{2}$. In any case, we have $A(x;\mu)\ge \alpha(\mu) \ge \mu_-:=\min\{ 1, \mu_i:\, i=1,...,p, \mu\in\cP\}$, i.e., in our case $\alpha(\mu)=\min\{\mu_1,\mu_2,1\}$ and $\mu_-=\frac18>0$.

This example allows us to check the quantitative behavior of our approach for a PPDE and also to compare different formulations (weak, ultra-weak).

\subsubsection{Standard weak form}\label{Sec:Num2UWF}
The standard weak form is also well-known to be induced by the bilinear form $b(w,v;\mu):=(A(\mu) w', v')_{L_2(0,1)} + (w, v)_{L_2(0,1)}$ which is coercive even on $\X=H^1(0,1)$, due to the reaction term. The coercivity constant is given by $\alpha(\mu)=\min\{\mu_1,\mu_2,1\}>\frac18$ as given above.

\subsubsection{Ultra-weak form}
The operator is easily seen to be self-adjoint, i.e., $\cB_\mu^{\circ,*} = \cB_\mu^{\circ}$. Moreover, $A(\mu)\in H^{1}_{\textrm{per}}(0,1)$ and $A(\mu)' \in \left(L_{\infty}(0,1)\right)$ for all $\mu \in\cP$. 
The bilinear form for the ultra-weak form is given by $b(w, v; \mu) = (w, \cB_\mu^{*}v)_{L_2(0,1)} =
-(w,  \left( A(\mu)  v' \right)')_{L_2(0,1)} + (w, v)_{L_2(0,1)}$ and $ \cB_\mu^{*}$ denotes the extension of $\cB_\mu^{\circ,*}$. 
The norm for the test space is parameter-dependent in this case, i.e.,  $\Vvert{v}\Vvert_\mu :=  \| \cB_\mu^{*}  v \|_{L_2(0,1)} = \| - \left( A(\mu) v' \right)' + v \|_{L_2(0,1)}$, which (for the one-dimensional case) can easily be shown to be equivalent to $ \| \cdot \|_{H^{2}(0,1)}$ for all $\mu \in\cP$, i.e.,
\begin{equation}\label{eq:NEqu}
	\alpha(\mu) \| v \|_{H^{2}(0,1)} 
	\le \Vvert{v}\Vvert_\mu \le C(\mu) \| v \|_{H^{2}(0,1)},
\end{equation}
where $C(\mu):= \max\{\| A(\mu)\|_{L_\infty(0,1)},1\}$. 
Hence, we obtain a parameter-independent test space $Y_{\mu} \equiv Y:= H_{\textrm{per}}^{2}(0,1)$. 
The norm equivalence \eqref{eq:NEqu} is crucial when replacing the dual norm $\Vert r_{\mu}^\theta \Vert_{(Y_{\mu})'}$ of the residual by $\Vert r_{\mu}^\theta \Vert_{H^{-2}(0,1)}$ and then using the wavelet norm equivalence \eqref{eq:weightedsum} for the latter term, i.e., $\sigma=-2$.

\subsubsection{Numerical comparisons}

\begin{figure}[htb]
	\centering
	\begin{tikzpicture}[scale = 0.95]
	\begin{semilogyaxis}[legend cell align={left}, grid=both,width=\textwidth,height=6cm,legend pos=south west,cycle list name=black white]
		\addplot[mark=*,red,line width=1pt] table [x=NoParam, y=H1Error] {figures/error_tc8_uw_w_new.dat};
		\addlegendentry{$H^1$-error weak}
		\addplot[mark=x,blue,line width=1pt] table [x=NoParam, y=EstH1w] {figures/error_tc8_uw_w_new.dat};
		\addlegendentry{estimate}
	\end{semilogyaxis}
	\end{tikzpicture} 
	\begin{tikzpicture}[scale = 0.95]
	\begin{semilogyaxis}[legend cell align={left}, grid=both,width=\textwidth,height=6cm,legend pos=north east,cycle list name=black white]
		\addplot[mark=*,red,line width=1pt] table [x=NoParam, y=L2Erroruw] {figures/error_tc8_uw_w_new.dat};
		\addlegendentry{$L_2$-error ultraweak}
		\addplot[mark=x,blue,line width=1pt] table [x=NoParam, y=Estuw] {figures/error_tc8_uw_w_new.dat};
		\addlegendentry{estimate}
	\end{semilogyaxis}
	\end{tikzpicture} 
	\caption{\label{fig:param_dependence}Error and error estimator for various parameters. Top: weak formulation, error measured in $H^1$; Bottom: ultra-weak formulation, error measured in $L_2$; The horizontal axis corresponds to the number of data points for $\mu$.}
\end{figure}
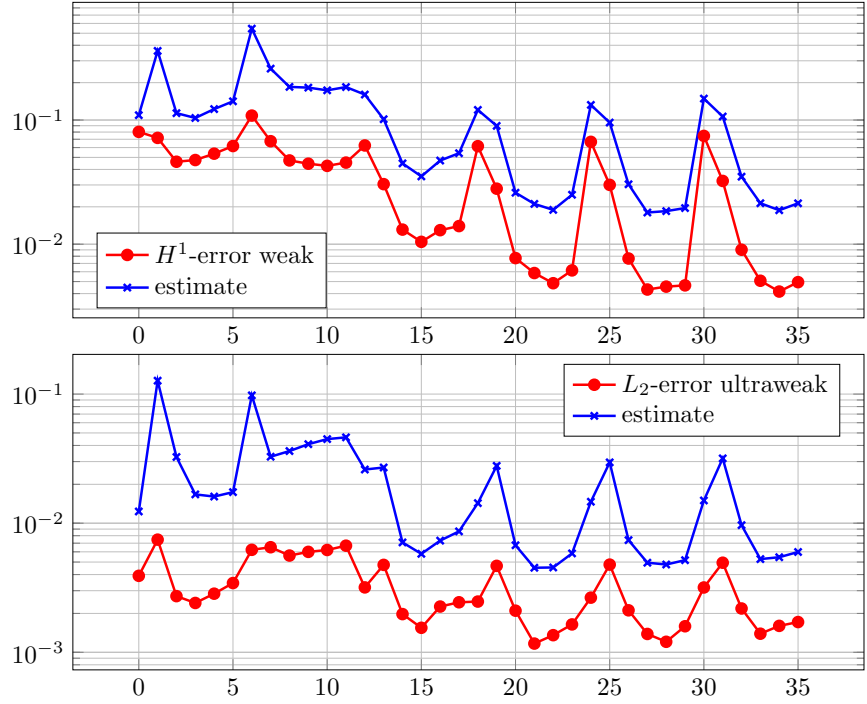

We report some results concerning the comparison of weak and ultra-weak formulations in Figure  \ref{fig:param_dependence}, where we show the exact error (red) and the wavelet-based error estimator (blue) for the standard variational form (top) and the ultra-weak form (bottom) for various parameters (consequently enumerated for the horizontal axis). In both cases, the wavelet-based error estimator is a strict upper bound. In the standard case, the effectivity is very good, which is due to the fact that the parameter-dependent coercivity constant is known and used here. For the ultra-weak form, the effectivity is slightly worse. We suppose that this is due to the replacement of the exact dual norm by $\Vert \cdot \Vert_{H^{-2}(0,1)}$. However, in both cases, we see that the error estimator follows the behavior of the exact error (the slopes of the curves are comparable).

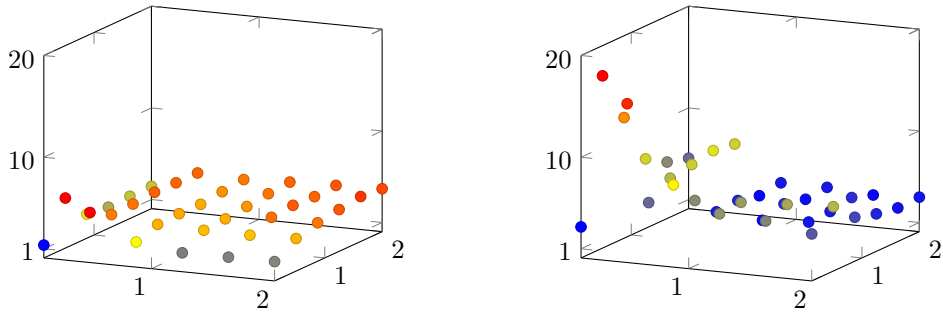
\begin{figure}[htb]
\begin{center}
\begin{tikzpicture}
	\begin{axis}[view = {25}{15},width=0.48\textwidth,zmin=0.1,zmax=20, ztick={1, 10, 20}]
	\addplot3+ [
	    domain=0.1:2,
	    domain y = 0.1:2,
	    samples = 6,
	    samples y = 6,
	    only marks,
	    scatter
	    ] table [x=mu1, y=mu2, z expr=\thisrow{EstH1w}/\thisrow{H1Error}] {figures/error_tc8_uw_w_new.dat};
	\end{axis}
\end{tikzpicture}
\hfill
\begin{tikzpicture}
	\begin{axis}[view = {25}{15},width=0.48\textwidth,zmin=0.1, zmax=20, ztick={1, 10, 20}]
	\addplot3+ [
	    domain=0.1:2,
	    domain y = 0.1:2,
	    samples = 6,
	    samples y = 6,
	    only marks,
	    scatter
	    ] table [x=mu1, y=mu2, z expr=\thisrow{Estuw}/\thisrow{L2Erroruw}] {figures/error_tc8_uw_w_new.dat};
	\end{axis}
\end{tikzpicture}
	\caption{\label{Fig_Effectivity}Effectivity of error estimator as scatter plot over parameters. Left: weak, $H^1$-norm; right: ultraweak, $L_2$-norm.}
\end{center}
\end{figure}

In order to show the effectivity in more detail, we show scatter plots of the effectivity in Figure \ref{Fig_Effectivity} for both cases (left: weak, right: ultra-weak). While in Figure \ref{fig:param_dependence} we just enumerated the points in parameter space, we show the data in Figure \ref{Fig_Effectivity} over the choice of $\mu\in\cP\subset\R^2$. We see that the effectivity is also good for many parameters with some outliers.

The above results also show that the wavelet-based PINN can be used for (nonlinear) model reduction as the online computation for a given parameter $\mu$ is just an evaluation of the PINN. In the above case, it is well-known that linear projection-based model reduction by the reduced basis method works extremely well, \cite{buffa-maday-patera-prudhomme-turinici-2012}. Hence, for those cases (allowing for a fast decay of the Kolmogorov $N$-width), there is no need for nonlinear model reduction. The situation changes, however, if one considers transport or wave-type problems. In that sense, the above results open the door towards certified nonlinear model reduction by wavelet-based PINNs for problems with slow decay of the Kolmogorov $N$-width. We will investigate this in future research.

\section{Summary and Outlook}
\label{sec:summary}
In this paper, we have used the well-known norm equivalence for scales of biorthogonal wavelets to approximate the dual norm of the residual of PPDEs up to any desired accuracy. This has been utilized within the framework of PINNs both for the definition of the loss function and for constructing an a posteriori error bound. This has been used for the standard variational formulation of elliptic (P)PDEs as well as for the ultra-weak formulation in $L_2$. Numerical experiments show very good performance of the wavelet-based PINN (with both variational formulations) in comparison even with a NN for approximating a known exact solution. We also see good effectivity (ratio of error bound and exact error) as long as the involved constants are approximated well.

This research shows that wavelet-based PINNs can efficiently be used for a certified approximation of PPDEs. This opens the door for additional research towards PINN-based nonlinear model reduction for further classes of PPDEs, in particular for those which are known to be problematic for linear model reduction.